\documentclass[12pt,leqno]{article}
\usepackage{amssymb}
\usepackage[mathscr]{eucal}
\usepackage{amsmath,amssymb,latexsym,theorem,bbm}
\usepackage{color}
\usepackage{appendix}

\newcommand{\SC}{\scriptstyle}

\newcommand{\CC}{\mathbb{C}}

\newcommand{\NN}{\mathbb{N}}
\newcommand{\RR}{\mathbb{R}}

\newcommand{\ZZ}{\mathbb{Z}}

\newcommand{\bA}{{\boldsymbol{A}}}

\newcommand{\bI}{{\boldsymbol{I}}}

\newcommand{\tm}{\widetilde{m}}

\newcommand{\bPhi}{{\boldsymbol{\Phi}}}

\newcommand{\bPsi}{{\boldsymbol{\Psi}}}
\newcommand{\bSigma}{{\boldsymbol{\Sigma}}}
\newcommand{\bzero}{{\boldsymbol{0}}}

\newcommand{\cB}{{\mathcal B}}

\newcommand{\cD}{{\mathcal D}}
\newcommand{\cE}{{\mathcal E}}

\newcommand{\cN}{{\mathcal N}}

\newcommand{\cX}{{\mathcal X}}

\newcommand{\cW}{{\mathcal W}}

\newcommand{\bcX}{\boldsymbol{\cX}}

\newcommand{\dd}{\mathrm{d}}
\newcommand{\ee}{\mathrm{e}}

\newcommand{\ii}{\mathrm{i}}

\newcommand{\EE}{\operatorname{\mathbb{E}}}
\newcommand{\PP}{{\operatorname{\mathbb{P}}}}
\newcommand{\QQ}{{\operatorname{\mathbb{Q}}}}

\newcommand{\var}{\operatorname{Var}}

\renewcommand{\Re}{{\operatorname{Re}}}
\renewcommand{\Im}{{\operatorname{Im}}}
\newcommand{\Res}{\operatorname{Res}}

\newcommand{\hh}{\widehat{h}}

\newcommand{\halpha}{\widehat{\alpha}}
\newcommand{\hbeta}{\widehat{\beta}}

\newcommand{\hvartheta}{\widehat{\vartheta}}

\newcommand{\talpha}{{\widetilde{\alpha}}}

\newcommand{\vare}{\varepsilon}

\renewcommand{\geq}{\geqslant}

\newcommand{\distr}{\stackrel{\cD}{\longrightarrow}}

\newcommand{\bbone}{\mathbbm{1}}

\newcommand{\proofend}{\hfill\mbox{$\Box$}}

\numberwithin{equation}{section}

\theoremstyle{change} \theorembodyfont{\em}
\newtheorem{Lem}{Lemma.}[section]
\newtheorem{Thm}[Lem]{Theorem.}

\newtheorem{Cor}[Lem]{Corollary.}

\theorembodyfont{\rm}

\begin{document}

\begin{center}
 {\bfseries\Large Nearly unstable family of stochastic processes given by} \\[2mm]
 {\bfseries\Large stochastic differential equations with time delay} \\[5mm]
 {\sc\large J\'anos Marcell $\text{Benke}^*$ \ and \ Gyula Pap}
\end{center}

\vskip0.2cm

\noindent
 Bolyai Institute, University of Szeged,
 Aradi v\'ertan\'uk tere 1, H--6720 Szeged, Hungary.

\noindent e--mail: jbenke@math.u-szeged.hu (J. M. Benke)

\noindent * Corresponding author.

\begin{center}
\textit{to Gyula}
\end{center}

\renewcommand{\thefootnote}{}
 \footnote{This research was supported by the EU-funded Hungarian grant EFOP-3.6.1-16-2016-00008.
The authors are supported by the Ministry of Human Capacities, Hungary
grant 20391-3/2018/FEKUSTRAT.}

\vspace*{-7mm}


\begin{abstract}
Let \ $a$ \ be a finite signed measure on \ $[-r, 0]$ \ with \ $r \in (0, \infty)$.
\ Consider a stochastic process \ $(X^{(\vartheta)}(t))_{t\in[-r,\infty)}$ \ given by a
 linear stochastic delay differential equation
 \[
   \dd X^{(\vartheta)}(t)
   = \vartheta \int_{[-r,0]} X^{(\vartheta)}(t + u) \, a(\dd u) \, \dd t + \dd W(t) ,
   \qquad t \in \RR_+ , 
 \]
 where \ $\vartheta \in \RR$ \ is a parameter and \ $(W(t))_{t\in\RR_+}$ \ is a
 standard Wiener process.
Consider a point \ $\vartheta \in \RR$, \ where this model is unstable in the sense that it is locally asymptotically Brownian functional with certain scalings \ $(r_{\vartheta,T})_{T\in(0,\infty)}$ \ satisfying \ $r_{\vartheta,T} \to 0$ \ as
 \ $T \to \infty$.
\ A family \ $\{(X^{(\vartheta_T)}(t))_{t\in[-r,T]} : T \in (0, \infty)\}$ \ is said to
 be nearly unstable as \ $T \to \infty$ \ if \ $\vartheta_T \to \vartheta$ \ as
 \ $T \to \infty$.
\ For every \ $\alpha \in \RR$, \ we prove convergence of the likelihood ratio
 processes of the nearly unstable families
 \ $\{(X^{(\vartheta+\alpha\hspace*{0.2mm}r_{\vartheta,T})}(t))_{t\in[-r,T]}
      : T \in (0, \infty)\}$
 \ as \ $T \to \infty$.
\ As a consequence, we obtain weak convergence of the maximum likelihood estimator
 \ $\halpha_T$ \ of \ $\alpha$ \ based on the observations
 \ $(X^{(\vartheta+\alpha\hspace*{0.2mm}r_{\vartheta,T})}(t))_{t\in[-r,T]}$ \ as
 \ $T \to \infty$.
\ It turns out that the limit distribution of \ $\halpha_T$ \ as \ $T \to \infty$ \ can
 be represented as the maximum likelihood estimator of a parameter of a process
 satisfying a stochastic differential equation without time delay.
\end{abstract}

\section{Introduction}

Research under the umbrellas of unstable, nearly unstable, unit root, near unit root, nonstationary, nearly nonstationary, integrated and near-integrated time series processes has received considerable attention in both the statistics and the econometric literature during the last 30 years, see, e.g., Chan \cite{Chan}. 

The aim of this paper is to show a phenomenon for certain nearly unstable family of stochastic processes given by stochastic differential equations with time delay which is known for nearly unstable sequences of AR(1) processes, see Bobkoski \cite{Bob}, Phillips \cite{Phi} and Chan and Wei \cite{ChanWei,ChanWei2}. 
Let us consider an AR(1) process
 \begin{equation}\label{AR1}
  \begin{cases}
   X^{(\beta)}_k = \beta X^{(\beta)}_{k-1} + \vare_k , \qquad
    k \in \NN := \{1, 2, \dots\} , \\
   X^{(\beta)}_0 = 0 ,
  \end{cases}
 \end{equation}
 where \ $\beta \in \RR$ \ is a parameter and \ $\vare_k$, \ $k \in \NN$,
 \ are independent standard normally distributed random variables.
Based on the asymptotic behavior of \ $\var(X^{(\beta)}_k) = 1 + \beta^2 + \cdots + \beta^{2(k-1)}$ \ as \ $k \to \infty$, \ the process \ $(X^{(\beta)})_{k\in\NN}$
\ is called stable, unstable or explosive if \ $|\beta| < 1$, \ $|\beta| = 1$ \ or \ $|\beta| > 1$, \ respectively.
For each \ $n \in \NN$, \ the least squares estimator (LSE) \ $\hbeta_n$ \ of \ $\beta$ \ based on a sample \ $\{X_0^{(\beta)}, \ldots, X_n^{(\beta)}\}$ \ is
 \[
   \hbeta_n
   = \frac{\sum_{k=1}^n X_{k-1}^{(\beta)} X_k^{(\beta)}}
          {\sum_{k=1}^n \bigl(X_{k-1}^{(\beta)}\bigr)^2} .
 \]
In the case of \ $|\beta| < 1$, \ the sequence \ $(\hbeta_n)_{n\in\NN}$ \ is asymptotically normal, namely, 
 \[
   \sqrt{n} (\hbeta_n - \beta) \distr \cN(0, 1 - \beta^2) \qquad \text{as \ $n \to \infty$,}
 \]
 see Mann and Wald \cite{ManWal} and Anderson \cite{And}.
In the case of \ $|\beta| = 1$, \ we have
 \[
   n (\hbeta_n - \beta) \distr \frac{\int_0^1 W(t) \, \dd W(t)}{\int_0^1 W(t)^2 \, \dd t} \qquad \text{as \ $n \to \infty$,}
 \]
 where \ $(W(t))_{t\in[0,1]}$ \ is a standard Wiener process, see White \cite{White}, Phillips \cite{Phi1, Phi2}, Greenwood and Wefelmeyer \cite[Example]{GreWef}.
In the case of \ $|\beta| > 1$, \ we have
 \[
   \beta^n (\hbeta_n - \beta) \distr \mathrm{Cauchy}(0, \beta^2 -1) \qquad \text{as \ $n \to \infty$,}
 \]
 see Anderson \cite{And} and White \cite{White}.
Moreover, considering for each \ $n \in \NN$ \ the statistical experiment \ $\cE_n$ \ connected with the observations
 \ $\{(X^{(\beta)}_0, \ldots, X^{(\beta)}_n) : \beta \in \RR\}$, \ the sequence \ $(\cE_n)_{n\in\NN}$ \ is locally asymptotically normal (LAN) if \ $|\beta| < 1$, \ it is locally asymptotically Brownian functional (LABF) if \ $\beta = 1$, \ and
 it is locally asymptotically mixed normal (LAMN) if \ $|\beta| > 1$, \ see, e.g., Greenwood and Wefelmeyer \cite[Example]{GreWef}. 

A sequence \ $(X^{(\beta^{(n)})})_{k\in\NN}$, \ $n \in \NN$, \ of AR(1) processes is called nearly unstable if \ $\beta^{(n)} \to \beta$ \ as \ $n \to \infty$, \ where \ $|\beta| = 1$.
\ In case of \ $\beta = 1$ \ and \ $\beta^{(n)} = 1 + \frac{h}{n}$ \ with some \ $h \in \RR$, \ we have
 \[
   n (\hbeta_n^{(n)} - \beta^{(n)}) \distr \frac{\int_0^1 Y(t) \, \dd W(t)}{\int_0^1 Y(t)^2 \, \dd t} \qquad \text{as \ $n \to \infty$,}
 \]
 where \ $(Y^{(h)}(t))_{t\in[0,1]}$ \ is a is a continuous time AR(1) process, i.e., an Ornstein--Uhlenbeck process, defined as a unique strong solution of the stochastic differential equation (SDE)
 \begin{equation}\label{OU}
  \begin{cases}
   \dd Y^{(h)}(t) = h Y^{(h)}(t) \dd t + \dd W(t), \qquad t \in [0, 1] , \\
   Y^{(h)}(0) = 0 ,
  \end{cases}
 \end{equation}
 see Bobkoski \cite{Bob}, Phillips \cite{Phi} and Chan and Wei \cite{ChanWei,ChanWei2}.
If we consider now \ $h$ \ as a parameter instead of \ $\beta^{(n)}$, \ then the LSE of \ $h = n (\hbeta^{(n)} - 1)$ \ is
 \[
   \hh_n = n (\hbeta_n^{(n)} - 1) = n (\hbeta_n^{(n)} - \beta^{(n)}) + h ,
 \]
 and we have
 \begin{equation} \label{LSEh}
  \hh_n \distr \frac{\int_0^1 Y^{(h)}(t) \, \dd W(t)}{\int_0^1 Y^{(h)}(t)^2 \, \dd t} + h
   = \frac{\int_0^1 Y^{(h)}(t) \, \dd Y^{(h)}(t)}{\int_0^1 Y^{(h)}(t)^2 \, \dd t} ,
 \end{equation}
 where the limit distribution in \eqref{LSEh} turns out to be the maximum likelihood estimator (MLE) of the parameter
\ $h$ \ in the model \eqref{OU} based on a sample \ $(Y^{(h)}(t))_{t\in[0,1]}$ \ (see, e.g., Arat\'o \cite{Arato}, or van der Meer, Pap and
 van Zuijlen \cite{MeerPapZui})). 
The same phenomenon appears in the nearly unstable AR($p$) model, see Jeganathan \cite{Jeg}, van der Meer,
 Pap and van Zuijlen \cite{MeerPapZui}, and Buchmann and Chan \cite{BucCha}.

In the present paper we consider a nearly unstable model for linear SDE with time
 delay.
For a linear SDE with time delay, the analog definition of unstability would be that the
 characteristic function has no roots with positive real part but has at least one root on
 the imaginary axis. 
However, investigating the local asymptotic properties of the likelihood function of this process,
 it turns out that the proper definition is slightly different, see Section 2 and
 Benke and Pap \cite{BenPap2017a}.
The main result of this paper is the convergence of the likelihood ratio processes of this nearly unstable model for linear SDE with time delay.
Based on this result, we show that in this nearly unstable model the same phenomenon appears which is described in this introduction earlier, namely, the limit distribution of the MLE of the model can
 be represented as the MLE of a parameter of a process
 satisfying a stochastic differential equation without time delay, which is in fact a multidimensional Ornstein--Uhlenbeck process.

\section{Preliminaries}

Let \ $\NN$, \ $\ZZ_+$, \ $\RR$, \ $\RR_+$, \ $\RR_{++}$, \ $\RR_{--}$ \ and \ $\CC$
 \ denote the sets of positive integers, non-negative integers, real numbers,
 non-negative real numbers, positive real numbers, negative real numbers and complex
 numbers, respectively.
Consider a linear stochastic delay differential equation (SDDE)
 \begin{equation}\label{SDDE}
  \begin{cases}
   \dd X^{(\vartheta)}(t)
   = \vartheta \int_{[-r,0]} X^{(\vartheta)}(t + u) \, a(\dd u) \, \dd t + \dd W(t) ,
    & t \in \RR_+ , \\
   X^{(\vartheta)}(t) = X_0(t) , & t \in [-r, 0] , 
  \end{cases}
 \end{equation}
 where \ $\vartheta \in \RR$ \ is a parameter, \ $r \in \RR_{++}$, \ $a$ \ is a finite
 signed measure on \ $[-r, 0]$ \ with \ $a \ne 0$, \ $(W(t))_{t\in\RR_+}$ \ is a
 standard Wiener process, and \ $(X_0(t))_{t\in[-r,0]}$ \ is a continuous stochastic
 process independent of \ $(W(t))_{t\in\RR_+}$.
\ The SDDE \eqref{SDDE} has a pathwise unique strong solution.
 
In the following, we recall some results from Benke and Pap \cite{BenPap2017a},
 where local asymptotic properties of the likelihood function has been studied.  
For all \ $T \in \RR_{++}$, \ let \ $\PP_{\vartheta,T}$ \ be the probability measure
 induced by \ $(X^{(\vartheta)}(t))_{t\in[-r,T]}$ \ on
 \ $(C([-r, T]), \cB(C([-r, T])))$. 
\ The Radon--Nikodym derivatives
 \ $\frac{\dd \PP_{\theta,T}}{\dd \PP_{\vartheta,T}}$, \ $\theta, \vartheta \in \RR$,
 \ can be derived from formula (7.139) in Section 7.6.4 of Liptser and Shiryaev
 \cite{LipShiI}.

\begin{Lem}\label{RN}
Let \ $\theta, \vartheta \in \RR$.
\ Then for all \ $T \in \RR_+$, \ the measures \ $\PP_{\theta,T}$ \ and
 \ $\PP_{\vartheta,T}$ \ are absolutely continuous with respect to each other, and
 \begin{align*}
  \log\left(\frac{\dd \PP_{\theta,T}}{\dd \PP_{\vartheta,T}}
            (X^{(\vartheta)}|_{[-r,T]})\right)
  &= (\theta - \vartheta) \int_0^T Y^{(\vartheta)}(t) \, \dd X^{(\vartheta)}(t)
     - \frac{1}{2} (\theta^2 - \vartheta^2) \int_0^T Y^{(\vartheta)}(t)^2 \, \dd t \\
  &= (\theta - \vartheta) \int_0^T Y^{(\vartheta)}(t) \, \dd W(t)
     - \frac{1}{2} (\theta - \vartheta)^2 \int_0^T Y^{(\vartheta)}(t)^2 \, \dd t
 \end{align*}
 with
 \[
   Y^{(\vartheta)}(t) := \int_{[-r,0]} X^{(\vartheta)}(t+u) \, a(\dd u) ,
   \qquad t \in \RR_+ .
 \]
Moreover, the process
 \ $\Bigl(\frac{\dd \PP_{\theta,T}}{\dd \PP_{\vartheta,T}}
          (X^{(\vartheta)}|_{[-r,T]})\Bigr)_{T\in\RR_+}$
 \ is a martingale.
\end{Lem}

The martingale property of the process
 \ $\Bigl(\frac{\dd \PP_{\theta,T}}{\dd \PP_{\vartheta,T}}
          (X^{(\vartheta)}|_{[-r,T]})\Bigr)_{T\in\RR_+}$
 \ is a consequence of Theorem 3.4 in Chapter III of Jacod and Shiryaev \cite{JSh}. 
We have the following simple corollary.

\begin{Cor}\label{RN_Cor}
For each \ $\vartheta \in \RR$, \ $T \in \RR_+$, \ $r_{\vartheta,T} \in \RR$ \ and
 \ $h_T \in \RR$, \ we have
 \[
   \log\left(\frac{\dd \PP_{\vartheta+r_{\vartheta,T}h_T,T}}{\dd \PP_{\vartheta,T}}
             (X^{(\vartheta)}|_{[-r,T]})\right)
   = h_T \Delta_{\vartheta,T} - \frac{1}{2} h_T^2 J_{\vartheta,T} ,
 \]
 with
 \[
   \Delta_{\vartheta,T}
   := r_{\vartheta,T} \int_0^T Y^{(\vartheta)}(t) \, \dd W(t) , \qquad
   J_{\vartheta,T}
   := r_{\vartheta,T}^2 \int_0^T Y^{(\vartheta)}(t)^2 \, \dd t .
 \]
\end{Cor}

The asymptotic behaviour of the solution of \eqref{SDDE} is connected with the
 so-called characteristic function \ $h_\vartheta : \CC \to \CC$, \ given by
 \begin{equation}\label{char1}
  h_\vartheta(\lambda)
  := \lambda - \vartheta \int_{[-r,0]} \ee^{\lambda u} \, a(\dd u) ,
  \qquad \lambda \in \CC ,
 \end{equation}
 and the set \ $\Lambda_\vartheta$ \ of the (complex) solutions of the so-called
 characteristic equation,
 \begin{equation}\label{char2}
  \lambda - \vartheta \int_{[-r,0]} \ee^{\lambda u} \, a(\dd u) = 0 .
 \end{equation}
For each \ $\lambda \in \Lambda_\vartheta$, \ denote by \ $\tm_\vartheta(\lambda)$
 \ the degree of the complex-valued polynomial
 \[
   P_{\vartheta,\lambda}(t)
   := \sum_{\ell=0}^{m_\vartheta(\lambda)-1}
       \frac{c_{\vartheta,\lambda,\ell}}{\ell!} t^\ell
 \]
 with
 \[
   c_{\vartheta,\lambda,\ell}
   := \int_{[-r, 0]}
       \underset{\SC z=\lambda}{\Res}
        \biggl(\frac{(z-\lambda)^\ell \ee^{zu}}{h_\vartheta(z)}\biggr)
       a(\dd u) ,
 \]
 where the degree of the zero polynomial is defined to be \ $-\infty$.
\ Put
 \begin{equation}\label{v*m*}
  v_\vartheta^*
  := \sup\{\Re(\lambda) : \lambda \in \Lambda_\vartheta , \;
                          \tm_\vartheta(\lambda) \geq 0\} , \qquad
  m_\vartheta^*
  := \max\{\tm_\vartheta(\lambda) : \lambda \in \Lambda_\vartheta, \;
                                    \Re(\lambda) = v_\vartheta^*\} ,
 \end{equation}
 where \ $\sup \emptyset := -\infty$ \ and \ $\max \emptyset := -\infty$.

The following theorem gives a sufficient condition for the LABF property of the family
 \begin{equation}\label{cET}
  (\cE_T)_{T\in\RR_{++}}
  := \big(C(\RR_+), \cB(C(\RR_+)),
          \{\PP_{\vartheta,T} : \vartheta \in \RR\}\big)_{T\in\RR_{++}}
 \end{equation}
 of statistical experiments.

\begin{Thm}\label{LABF}
If \ $\vartheta \in \RR$ \ with \ $v_\vartheta^* = 0$, \ then the family
 \ $(\cE_T)_{T\in\RR_{++}}$ \ of statistical experiments given in \eqref{cET} is LABF at
 \ $\vartheta$, \ namely,
 \begin{equation}\label{LAQD}
  (\Delta_{\vartheta,T}, J_{\vartheta,T}) \distr (\Delta_\vartheta, J_\vartheta)
  \qquad \text{as \ $T \to \infty$}
 \end{equation}
 with scaling \ $r_{\vartheta,T} = T^{-m_\vartheta^*-1}$ \ and with
 \begin{gather*}
  \Delta_\vartheta
  = \sum_{\underset{\SC\tm_\vartheta(\lambda)=m_\vartheta^*}
                   {\lambda \in \Lambda_\vartheta\cap(\ii\RR)}}
     c_{\vartheta,\lambda,m_\vartheta^*}
     \int_0^1
      \cW_{\Im(\lambda),m_\vartheta^*}(s)
      \, \dd \overline{\cW_{\Im(\lambda)}(s)} , \\
  J_\vartheta
  = \sum_{\underset{\SC\tm_\vartheta(\lambda)=m_\vartheta^*}
                   {\lambda \in \Lambda_\vartheta\cap(\ii\RR)}}
     |c_{\vartheta,\lambda,m_\vartheta^*}|^2
     \int_0^1 |\cW_{\Im(\lambda),m_\vartheta^*}(s)|^2 \, \dd s ,
 \end{gather*}
 with
 \[
   \cW_\varphi
   := \begin{cases}
       \cW_0, & \text{if \ $\varphi = 0$,} \\
       \frac{1}{\sqrt{2}} \bigl(\cW_\varphi^\Re + \ii \cW_\varphi^\Im\bigr),
        & \text{if \ $\varphi \in \RR_{++}$,} \\[1mm]
       \overline{\cW_{-\varphi}},
        & \text{if \ $\varphi \in \RR_{--}$,}
      \end{cases}
 \]
 where \ $\{\cW_0, \cW_\varphi^\Re, \cW_\varphi^\Im : \varphi \in \RR_{++}\}$ \ are
 independent standard Wiener processes, and
 \[
   \cW_{\varphi,\ell}(s)
   := \frac{1}{\ell!} \int_0^s (s - u)^\ell \, \dd \cW_\varphi(u) , \qquad
   s \in \RR_+ , \quad \varphi \in \RR , \quad \ell \in \ZZ_+ .
 \]
Particularly, if \ $a([-r, 0]) \ne 0$, \ then \ $v_0^* = 0$, \ $m_0^* = 0$, \ and the
 family \ $(\cE_T)_{T\in\RR_{++}}$ \ of statistical experiments given in \eqref{cET} is
 LABF at \ $0$ \ with scaling \ $r_{0,T} = T^{-1}$, \ $T \in \RR_{++}$, \ and with 
 \[
   \Delta_0 = a([-r, 0]) \int_0^1 \cW_0(s) \, \dd \cW_0(s) , \qquad
   J_0 = a([-r, 0])^2 \int_0^1 \cW_0(s)^2 \, \dd s . 
 \]
\end{Thm}

Note that
 \ $c_{\vartheta,\overline{\lambda},m_\vartheta^*}
    = \overline{c_{\vartheta,\lambda,m_\vartheta^*}}$
 \ for all \ $\lambda \in \Lambda_\vartheta$, \ hence, if
 \ $0 \in \Lambda_\vartheta$, \ then \ $c_{\vartheta,0,m_\vartheta^*} \in \RR$. 
\ Moreover,
 \begin{equation}\label{conj}
  c_{\vartheta,\overline{\lambda},m_\vartheta^*}
    \int_0^1
     \cW_{\Im(\overline{\lambda}),m_\vartheta^*}(s)
     \, \dd \overline{\cW_{\Im(\overline{\lambda})}(s)}
  = \overline{c_{\vartheta,\lambda,m_\vartheta^*}}
              \int_0^1
               \overline{\cW_{\Im(\lambda),m_\vartheta^*}(s)}
               \, \dd \cW_{\Im(\lambda)}(s)
 \end{equation}
 for all \ $\lambda \in \Lambda_\vartheta$, \ hence \ $\Delta_\vartheta$ \ is almost
 surely real-valued. 
 
Based on Theorem \ref{LABF}, we say that the process given by \eqref{SDDE} is unstable if \ $v_\vartheta^* = 0$.
\ Note that the family \ $(\cE_T)_{T\in\RR_{++}}$ \ of statistical experiments given in \eqref{cET} is LAN if \ $v_\vartheta^* < 0$, \ and, under some additional conditions, LAMN or periodically locally asymptotic mixed normal (PLAMN) if \ $v_\vartheta^* > 0$, \ see Benke and Pap \cite{BenPap2017a}.

\section{Nearly unstable models}
\label{main}

Let \ $\vartheta \in \RR$ \ with \ $v_\vartheta^* = 0$, \ and consider the scaling \ $r_{\vartheta,T} = T^{-m_\vartheta^*-1}$. 
\ By Theorem \ref{LABF}, for each \ $\alpha \in \RR$, \ the family
 \ $\{(X^{(\vartheta+\alpha\hspace*{0.2mm}r_{\vartheta,T})}(t))_{t\in[-r,T]}
      : T \in (0, \infty)\}$
 \ is nearly unstable \ as \ $T \to \infty$.
\ In order to describe the asymptotic behavior of the likelihood ratio process of
 \ $(X^{(\vartheta+\alpha\hspace*{0.2mm}r_{\vartheta,T})}(t))_{t\in[-r,T]}$ \ as
 \ $T \to \infty$, \ we need certain stochastic processes.
For each \ $\vartheta \in \RR$ \ and \ $\lambda \in \Lambda_\vartheta \cap (\ii \RR)$,
 \ consider the linear SDE (without delay)
 \begin{equation}\label{SDElc}
  \begin{cases}
   \dd \cX_{\vartheta,\lambda,0}^{(\alpha)}(t)
   = \alpha \hspace*{0.2mm} c_{\vartheta,\lambda,m^*_\vartheta}
     \hspace*{0.2mm}
     \cX_{\vartheta,\lambda,m^*_\vartheta}^{(\alpha)}(t) 
     \, \dd t
     + \dd \cW_{\Im(\lambda)}(t) ,
    & t \in [0, 1] , \\[1mm]
   \dd \cX_{\vartheta,\lambda,\ell}^{(\alpha)}(t)
   = \cX_{\vartheta,\lambda,\ell-1}^{(\alpha)}(t) \, \dd t , \qquad
    \ell \in \{1, \ldots, m^*_\vartheta\} ,
    & t \in [0, 1] , \\[1mm]
   \cX_{\vartheta,\lambda,\ell}^{(\alpha)}(0) = 0 , \qquad
    \ell \in \{0, 1, \ldots, m^*_\vartheta\} ,
  \end{cases}
 \end{equation}
 where \ $\alpha \in \RR$ \ is a parameter.
The SDE \eqref{SDElc} has a pathwise unique strong solution.
The processes \ $(\cX_{\vartheta,\lambda,\ell}^{(\alpha)}(t))_{t\in[0,1]}$,
 \ $\ell \in \{0, 1, \ldots, m^*_\vartheta\}$, \ are real-valued if
 \ $\Im(\lambda) = 0$, \ and they are complex-valued if \ $\Im(\lambda) \ne 0$.
\ Moreover, for each \ $\ell \in \{0, 1, \ldots, m^*_\vartheta\}$ \ and
 \ $t \in [0, 1]$, \ we have
 \ $\cX_{\vartheta,\overline{\lambda},\ell}^{(\alpha)}(t)
    = \overline{\cX_{\vartheta,\lambda,\ell}^{(\alpha)}(t)}$.
\ Put
 \[
   \bcX^{(\alpha)}_\vartheta(t)
   := \bigl(\bcX^{(\alpha)}_{\vartheta,\lambda}(t)
            : \lambda \in \Lambda_\vartheta \cap (\ii \RR_+),
              \tm_\vartheta(\lambda) = m^*_\vartheta\bigr) , \qquad t \in [0, 1] ,
 \]
 where
 \[
   \bcX^{(\alpha)}_{\vartheta,\lambda}(t)
   := \begin{cases}
       \bigl(\cX_{\vartheta,0,0}^{(\alpha)}(t), \ldots,
             \cX_{\vartheta,0,m^*_\vartheta}^{(\alpha)}(t)\bigr)^\top ,
        & \text{if \ $\Im(\lambda) = 0$,} \\
       \bigl(\bPhi(\cX_{\vartheta,\lambda,0}^{(\alpha)}(t))^\top, \ldots,
             \bPhi(\cX_{\vartheta,\lambda,m^*_\vartheta}^{(\alpha)}(t))^\top
       \bigr)^\top ,
        & \text{if \ $\Im(\lambda) \in \RR_{++}$,}
      \end{cases}
 \]
 with
 \[
   \bPhi(z) := \begin{bmatrix}
                 \Re(z) \\
                 \Im(z)
                \end{bmatrix} , \qquad z \in \CC .
 \]
If \ $\Im(\lambda) = 0$, \ then the \ $(m^*_\vartheta + 1)$-dimensional real-valued
 process \ $(\bcX^{(\alpha)}_{\vartheta,0}(t))_{t\in[0,1]}$ \ is the pathwise unique
 strong solution of the linear SDE
 \begin{equation}\label{SDE0}
  \dd \bcX^{(\alpha)}_{\vartheta,0}(t)
  = \bA^{(\alpha)}_{\vartheta,0} \hspace*{0.2mm} \bcX^{(\alpha)}_{\vartheta,0}(t)
    \, \dd t
    + \bSigma_0 \, \dd \cW_0(t) , \qquad t \in [0, 1] ,
 \end{equation}
 with
 \[
   \bA^{(\alpha)}_{\vartheta,0}
   := \begin{bmatrix}
       0 & 0 & \cdots & 0 & \alpha \hspace*{0.2mm} c_{\vartheta, 0, m^*_\vartheta} \\
       1 & 0 & \cdots & 0 & 0 \\
       0 & 1 & \cdots & 0 & 0 \\
       \vdots & \vdots & \ddots & \vdots & \vdots \\
       0 & 0 & \cdots & 1 & 0
      \end{bmatrix} , \qquad
   \bSigma_0
   := \begin{bmatrix}
       1 \\
       0 \\
       0 \\
       \vdots \\
       0
      \end{bmatrix} .
 \]
Let \ $\PP^{(\alpha)}_{\vartheta,0}$ \ be the probability measure induced by the
 process \ $(\bcX^{(\alpha)}_{\vartheta,0}(t))_{t\in[0,1]}$ \ on the space
 \ $(C([0, 1])^{m^*_\vartheta+1}, \cB(C([0, 1])^{m^*_\vartheta+1}))$. 
\ Applying again formula (7.139) in Section 7.6.4 of Liptser and Shiryaev
 \cite{LipShiI}, we obtain that for all \ $\alpha, \talpha \in \RR$, \ the measures
 \ $\PP^{(\alpha)}_{\vartheta,0}$ \ and \ $\PP^{(\talpha)}_{\vartheta,0}$ \ are
 absolutely continuous with respect to each other, and
 \begin{align*}
  \log\left(\frac{\dd\PP^{(\talpha)}_{\vartheta,0}}{\dd\PP^{(\alpha)}_{\vartheta,0}}
   (\bcX^{(\alpha)}_{\vartheta,0})\right)
  &= \int_0^1 
      \bigl((\bA^{(\talpha)}_{\vartheta,0} - \bA^{(\alpha)}_{\vartheta,0})
            \bcX^{(\alpha)}_{\vartheta,0}(t)\bigr)^\top
      (\bSigma_0 \bSigma_0^\top)^\ominus
      \, \dd \bcX^{(\alpha)}_{\vartheta,0}(t) \\
  &\quad
     - \frac{1}{2}
       \int_0^1 
        \bigl((\bA^{(\talpha)}_{\vartheta,0} - \bA^{(\alpha)}_{\vartheta,0})
              \bcX^{(\alpha)}_{\vartheta,0}(t)\bigr)^\top
        (\bSigma_0 \bSigma_0^\top)^\ominus
        (\bA^{(\talpha)}_{\vartheta,0} + \bA^{(\alpha)}_{\vartheta,0})
        \bcX^{(\alpha)}_{\vartheta,0}(t)
        \, \dd t ,
 \end{align*}
 where \ $(\bSigma_0 \bSigma_0^\top)^\ominus$ \ denotes the generalized inverse of
 \ $\bSigma_0 \bSigma_0^\top$. 
\ We have \ $(\bSigma_0 \bSigma_0^\top)^\ominus = \bSigma_0 \bSigma_0^\top$ \ and
 \begin{align*}
  \bigl((\bA^{(\talpha)}_{\vartheta,0} - \bA^{(\alpha)}_{\vartheta,0})
        \bcX^{(\alpha)}_{\vartheta,0}(t)\bigr)^\top
  \bSigma_0
  &= (\talpha - \alpha) c_{\vartheta,0,m^*_\vartheta}
     \cX^{(\alpha)}_{\vartheta,0,m^*_\vartheta}(t) , \\
  \bSigma_0^\top (\bA^{(\talpha)}_{\vartheta,0} + \bA^{(\alpha)}_{\vartheta,0})
  \bcX^{(\alpha)}_{\vartheta,0}(t)
  &= (\talpha + \alpha) c_{\vartheta,0,m^*_\vartheta} 
     \cX^{(\alpha)}_{\vartheta,0,m^*_\vartheta}(t), \\
  \bSigma_0^\top \, \dd \bcX^{(\alpha)}_{\vartheta,0}(t)
  &= \dd \cX^{(\alpha)}_{\vartheta,0,0}(t) ,
 \end{align*}
 hence, using the SDE \eqref{SDE0},
 \begin{align*}
  \log\left(\frac{\dd\PP^{(\talpha)}_{\vartheta,0}}{\dd\PP^{(\alpha)}_{\vartheta,0}}
   (\bcX^{(\alpha)}_{\vartheta,0})\right)
  &= (\talpha - \alpha) c_{\vartheta, 0, m^*_\vartheta}
      \int_0^1
       \cX^{(\alpha)}_{\vartheta,0,m^*_\vartheta}(t)
       \, \dd \cX^{(\alpha)}_{\vartheta,0,0}(t) \\
  &\quad
     - \frac{1}{2}
     (\talpha^2 - \alpha^2) c_{\vartheta,0,m^*_\vartheta}^2
      \int_0^1 \cX^{(\alpha)}_{\vartheta,0,m^*_\vartheta}(t)^2 \, \dd t \\
  &= (\talpha - \alpha) \Delta^{(\alpha)}_{\vartheta,0}
     - \frac{1}{2} (\talpha - \alpha)^2 J^{(\alpha)}_{\vartheta,0}
 \end{align*}
 with
 \[
   \Delta^{(\alpha)}_{\vartheta,0}
   := c_{\vartheta,0,m^*_\vartheta}
      \int_0^1 \cX^{(\alpha)}_{\vartheta,0,m^*_\vartheta}(t) \, \dd \cW_0(t) , \qquad
   J^{(\alpha)}_{\vartheta,0}
   := c_{\vartheta,0,m^*_\vartheta}^2
      \int_0^1 \cX^{(\alpha)}_{\vartheta,0,m^*_\vartheta}(t)^2 \, \dd t .
 \]
If \ $\Im(\lambda) \in \RR_{++}$, \ then the \ $2(m^*_\vartheta+1)$-dimensional
 real-valued process \ $(\bcX^{(\alpha)}_{\vartheta,\lambda}(t))_{t\in[0,1]}$ \ is the
 pathwise unique strong solution of the linear SDE
 \begin{equation}\label{SDEl}
   \dd \bcX^{(\alpha)}_{\vartheta,\lambda}(t)
   = \bA^{(\alpha)}_{\vartheta,\lambda} \hspace*{0.2mm}
     \bcX^{(\alpha)}_{\vartheta,\lambda}(t)
     \, \dd t
     + \frac{1}{\sqrt{2}}
       \bSigma_\lambda
       \begin{bmatrix}
        \dd \cW_{\Im(\lambda)}^\Re(t) \\
        \dd \cW_{\Im(\lambda)}^\Im(t)
       \end{bmatrix} , \qquad t \in [0, 1] ,
 \end{equation}
 with
 \[
   \bA^{(\alpha)}_{\vartheta,\lambda}
   := \begin{bmatrix}
       \bzero_{2\times2} & \cdots & \bzero_{2\times2}
        & \alpha \bPsi(c_{\vartheta,\lambda,m^*_\vartheta}) \\
       \bI_2 & \cdots & \bzero_{2\times2} & \bzero_{2\times2} \\
       \vdots & \ddots & \vdots & \vdots \\
       \bzero_{2\times2} & \cdots & \bI_2 & \bzero_{2\times2}
      \end{bmatrix} , \qquad
   \bSigma_\lambda
   := \begin{bmatrix}
       \bI_2 \\
       \bzero_{2\times2} \\
       \vdots \\
       \bzero_{2\times2}
      \end{bmatrix} ,
 \]
 where
 \[
   \bPsi(z) := \begin{bmatrix}
                \Re(z) & -\Im(z) \\
                \Im(z) & \Re(z)
               \end{bmatrix} , \quad z \in \CC , \qquad
   \bzero_{2\times2} := \begin{bmatrix}
                         0 & 0 \\
                         0 & 0
                        \end{bmatrix} , \qquad
   \bI_2 := \begin{bmatrix}
             1 & 0 \\
             0 & 1
            \end{bmatrix} .
 \]
Indeed, we have
 \begin{equation}\label{bPsibPhi}
  \bPsi(z_1) \bPhi(z_2)
  = \bPhi(z_1 z_2) , \qquad z_1, z_2 \in \CC ,
 \end{equation}
 and hence
 \ $\bPhi(c_{\vartheta,\lambda,m^*_\vartheta} \hspace*{0.2mm}
          \cX_{\vartheta,\lambda,m^*_\vartheta}^{(\alpha)}(t))
    = \bPsi(c_{\vartheta,\lambda,m^*_\vartheta})
      \bPhi(\cX_{\vartheta,\lambda,m^*_\vartheta}^{(\alpha)}(t))$.
      
Let \ $\PP^{(\alpha)}_{\vartheta,\lambda}$ \ be the probability measure induced by the
 process \ $(\bcX^{(\alpha)}_{\vartheta,\lambda}(t))_{t\in[0,1]}$ \ on the space
 \ $(C([0, 1])^{2(m^*_\vartheta+1)}, \cB(C([0, 1])^{2(m^*_\vartheta+1)}))$. 
\ Applying again formula (7.139) in Section 7.6.4 of Liptser and Shiryaev
 \cite{LipShiI}, we obtain that for all \ $\alpha, \talpha \in \RR$, \ the measures
 \ $\PP^{(\alpha)}_{\vartheta,\lambda}$ \ and \ $\PP^{(\talpha)}_{\vartheta,\lambda}$
 \ are absolutely continuous with respect to each other, and
 \begin{align*}
  \log\left(\frac{\dd\PP^{(\talpha)}_{\vartheta,\lambda}}
                 {\dd\PP^{(\alpha)}_{\vartheta,\lambda}}
            (\bcX^{(\alpha)}_{\vartheta,\lambda})\right)
  &= \int_0^1 
      \bigl((\bA^{(\talpha)}_{\vartheta,\lambda} - \bA^{(\alpha)}_{\vartheta,\lambda})
            \bcX^{(\alpha)}_{\vartheta,\lambda}(t)\bigr)^\top
      \left(\frac{1}{2} \bSigma_\lambda \bSigma_\lambda^\top\right)^\ominus
      \, \dd \bcX^{(\alpha)}_{\vartheta,\lambda}(t) \\
  &\quad
     - \frac{1}{2}
       \int_0^1 
        \bigl((\bA^{(\talpha)}_{\vartheta,\lambda}
               - \bA^{(\alpha)}_{\vartheta,\lambda})
              \bcX^{(\alpha)}_{\vartheta,\lambda}(t)\bigr)^\top
        \left(\frac{1}{2} \bSigma_\lambda \bSigma_\lambda^\top\right)^\ominus
        (\bA^{(\talpha)}_{\vartheta,\lambda} + \bA^{(\alpha)}_{\vartheta,\lambda})
        \bcX^{(\alpha)}_{\vartheta,\lambda}(t)
        \, \dd t ,
 \end{align*}
 where
 \ $\left(\frac{1}{2} \bSigma_\lambda \bSigma_\lambda^\top\right)^\ominus
    = 2 \bSigma_\lambda \bSigma_\lambda^\top$
 \ and
 \begin{align*}
  \bigl((\bA^{(\talpha)}_{\vartheta,\lambda} - \bA^{(\alpha)}_{\vartheta,\lambda})
        \bcX^{(\alpha)}_{\vartheta,\lambda}(t)\bigr)^\top
  \bSigma_\lambda
  &= (\talpha - \alpha)
     \begin{bmatrix}
      \Re(\cX^{(\alpha)}_{\vartheta,\lambda,m^*_\vartheta}(t)) \\
      \Im(\cX^{(\alpha)}_{\vartheta,\lambda,m^*_\vartheta}(t))
     \end{bmatrix}^\top
     \bPsi(c_{\vartheta,\lambda,m^*_\vartheta})^\top \\
  &= (\talpha - \alpha) \bPhi(\cX^{(\alpha)}_{\vartheta,\lambda,m^*_\vartheta}(t))^\top
     \bPsi(c_{\vartheta,\lambda,m^*_\vartheta})^\top , \\
  \bSigma_\lambda^\top 
  (\bA^{(\talpha)}_{\vartheta,\lambda} + \bA^{(\alpha)}_{\vartheta,\lambda})
  \bcX^{(\alpha)}_{\vartheta,\lambda}(t)
  &= (\talpha + \alpha) \bPsi(c_{\vartheta,\lambda,m^*_\vartheta})
     \begin{bmatrix}
      \Re(\cX^{(\alpha)}_{\vartheta,\lambda,m^*_\vartheta}(t)) \\
      \Im(\cX^{(\alpha)}_{\vartheta,\lambda,m^*_\vartheta}(t))
     \end{bmatrix} \\
  &= (\talpha + \alpha) \bPsi(c_{\vartheta,\lambda,m^*_\vartheta})
     \bPhi(\cX^{(\alpha)}_{\vartheta,\lambda,m^*_\vartheta}(t)) , \\
  \bSigma_\lambda^\top \, \dd \bcX^{(\alpha)}_{\vartheta,\lambda}(t)
  &= \begin{bmatrix}
      \dd\Re(\cX^{(\alpha)}_{\vartheta,\lambda,0}(t)) \\
      \dd\Im(\cX^{(\alpha)}_{\vartheta,\lambda,0}(t))
     \end{bmatrix}
   = \bPhi(\dd \cX^{(\alpha)}_{\vartheta,\lambda,0}(t)) ,
 \end{align*}
 hence, using the SDE \eqref{SDEl}, the identities \eqref{bPsibPhi} and
 \begin{gather*}
  \bPsi(z)^\top \bPsi(z) = |z|^2 \bI_2 , \qquad z \in \CC , \\
  \bPhi(z_1)^\top \bPhi(z_2) = \Re(z_1 \overline{z_2}) , \qquad z_1, z_2 \in \CC , \\
  \bPhi(z)^\top \bPhi(z) = |z|^2 , \qquad z \in \CC ,
 \end{gather*}
 we obtain 
 \begin{align*}
  \log\left(\frac{\dd\PP^{(\talpha)}_{\vartheta,\lambda}}
                 {\dd\PP^{(\alpha)}_{\vartheta,\lambda}}
            (\bcX^{(\alpha)}_{\vartheta,\lambda})\right)
  &= 2 (\talpha - \alpha)
     \int_0^1
      \bPhi(\cX^{(\alpha)}_{\vartheta,\lambda,m^*_\vartheta}(t))^\top
      \bPsi(c_{\vartheta,\lambda,m^*_\vartheta})^\top
      \bPhi(\dd \cX^{(\alpha)}_{\vartheta,\lambda,0}(t)) \\
  &\quad
     - (\talpha^2 - \alpha^2) |c_{\vartheta,0,m^*_\vartheta}|^2
       \int_0^1 |\cX^{(\alpha)}_{\vartheta,0,m^*_\vartheta}(t)|^2 \, \dd t \\
  &= (\talpha - \alpha) \Delta^{(\alpha)}_{\vartheta,\lambda}
     - \frac{1}{2} (\talpha - \alpha)^2 J^{(\alpha)}_{\vartheta,\lambda}
 \end{align*}
 with
 \[
   \Delta^{(\alpha)}_{\vartheta,\lambda}
   := 2 \Re\left(\int_0^1
                  c_{\vartheta,\lambda,m^*_\vartheta}
                  \cX^{(\alpha)}_{\vartheta,\lambda,m^*_\vartheta}(t)
                  \, \dd \overline{\cW_{\Im(\lambda)}(t)}\right) , \qquad
   J^{(\alpha)}_{\vartheta,\lambda}
   := 2 |c_{\vartheta,\lambda,m^*_\vartheta}|^2
      \int_0^1 |\cX^{(\alpha)}_{\lambda,m^*_\vartheta}(t)|^2 \, \dd t .
 \]
Let \ $\PP^{(\alpha)}_\vartheta$ \ be the probability measure induced by the process
 \ $(\bcX^{(\alpha)}_\vartheta(t))_{t\in[0,1]}$ \ on the space
 \ $(C([0, 1])^{d_\vartheta}, \cB(C([0, 1])^{d_\vartheta}))$ \ with
 \[
   d_\vartheta
   := \sum_{\underset{\SC\tm_\vartheta(\lambda)=m_\vartheta^*}
           {\lambda \in \Lambda_\vartheta\cap(\ii\RR_+)}}
       (m^*_\vartheta + 1) .
 \]
We obtain that for all \ $\alpha, \talpha \in \RR$, \ the measures \ $\PP^{(\alpha)}$
 \ and \ $\PP^{(\talpha)}$ \ are absolutely continuous with respect to each other, and
 \[
   \log\left(\frac{\dd\PP^{(\talpha)}}{\dd\PP^{(\alpha)}}(\bcX^{(\alpha)})\right)
   = \sum_{\underset{\SC\tm_\vartheta(\lambda)=m_\vartheta^*}
          {\lambda \in \Lambda_\vartheta\cap(\ii\RR_+)}}
      \log\left(\frac{\dd \QQ^{(\talpha)}_\lambda}{\dd \QQ^{(\alpha)}_\lambda}
                (\bcX^{(\alpha)}_\lambda)\right)
   = (\talpha - \alpha) \Delta^{(\alpha)}_\vartheta
     - \frac{1}{2} (\talpha - \alpha)^2 J^{(\alpha)}_\vartheta
 \]
 with
 \[
   \Delta^{(\alpha)}_\vartheta
   := \sum_{\underset{\SC\tm_\vartheta(\lambda)=m_\vartheta^*}
           {\lambda \in \Lambda_\vartheta\cap(\ii\RR_+)}}
       \Delta^{(\alpha)}_{\vartheta,\lambda} , \qquad
   J^{(\alpha)}_\vartheta
   := \sum_{\underset{\SC\tm_\vartheta(\lambda)=m_\vartheta^*}
           {\lambda \in \Lambda_\vartheta\cap(\ii\RR_+)}}
       J^{(\alpha)}_{\vartheta,\lambda} ,
 \]
 since independence of the Wiener processes
 \ $\{\cW_0, \cW_\varphi^\Re, \cW_\varphi^\Im : \varphi \in \RR_{++}\}$ \ yields
 independence of the processes
 \ $\{\bcX^{(\alpha)}_{\vartheta,\lambda} : \lambda \in \RR_+\}$.
 
\begin{Thm}\label{conv}
If \ $\vartheta \in \RR$ \ with \ $v_\vartheta^* = 0$, \ then for each
 \ $\alpha \in \RR$,
 \begin{equation}\label{convDJ}
   (\Delta_{\vartheta+\alpha\hspace*{0.1mm}r_{\vartheta,T},\,T},
    J_{\vartheta+\alpha\hspace*{0.1mm}r_{\vartheta,T},\,T})
   \distr
   \bigl(\Delta^{(\alpha)}_\vartheta, J^{(\alpha)}_\vartheta\bigr) \qquad
   \text{as \ $T \to \infty$.}
 \end{equation}
Consequently, the family
 \ $\bigl(C([-r,T]), \cB(C([-r,T])),
          \{\PP_{\vartheta+(\alpha+h)r_{\vartheta,T},\,T} 
            : h \in \RR\})_{T\in\RR_{++}}$
 \ of statistical experiments converge to the statistical experiment
 \ $(\RR^2, \cB(\RR^2), \{\QQ^{(\alpha)}_{\vartheta,h} : h \in \RR\})$ \ as
 \ $T \to \infty$, \ where
 \[
   \QQ^{(\alpha)}_{\vartheta,h}(B)
   := \EE\left(\exp\left\{h \Delta^{(\alpha)}_\vartheta
                          - \frac{1}{2} h^2 J^{(\alpha)}_\vartheta \right\}
               \bbone_B\bigl(\Delta^{(\alpha)}_\vartheta,
                             J^{(\alpha)}_\vartheta\bigr)\right) ,
   \qquad B \in \cB(\RR^2) ,
 \]
 in the sense that for each base \ $h_0 \in \RR$, \ the finite dimensional 
 distributions of the likelihood ratio process
 \ $\Bigl(\frac{\dd\PP_{\vartheta+(\alpha+h)r_{\vartheta,T},\,T}}
               {\dd\PP_{\vartheta+(\alpha+h_0)r_{\vartheta,T},\,T}}\Bigr)_{h\in\RR}$
 \ under \ $\PP_{\vartheta+(\alpha+h_0)r_{\vartheta,T},\,T}$ \ converge to the finite
 dimensional distributions of the likelihood ratio process
 \ $\Bigl(\frac{\dd\QQ^{(\alpha)}_{\vartheta,h}}
               {\dd\QQ^{(\alpha)}_{\vartheta,h_0}}\Bigr)_{h\in\RR}$
 \ under \ $\QQ^{(\alpha)}_{\vartheta,h_0}$ \ as \ $T \to \infty$.
\end{Thm}

Note that the probability measure \ $\QQ^{(\alpha)}_{\vartheta,0}$ \ is the 
 distribution of the random vector
 \ $\bigl(\Delta^{(\alpha)}_\vartheta, J^{(\alpha)}_\vartheta\bigr)$, \ which is
 concentrated on \ $\RR \times \RR_+$. 
\ Moreover, for each \ $h \in \RR$, \ the probability measures
 \ $\QQ^{(\alpha)}_{\vartheta,h}$ \ and \ $\QQ^{(\alpha)}_{\vartheta,0}$ \ are
 equivalent with
 \[
   \frac{\dd\QQ^{(\alpha)}_{\vartheta,h}}{\dd\QQ^{(\alpha)}_{\vartheta,0}}(\Delta, J)
   = \exp\left\{h \Delta - \frac{1}{2} h^2 J\right\} , \qquad (\Delta, J) \in \RR^2 .
 \] 
 
\noindent
\textbf{Proof of Theorem \ref{conv}.} 
By Corollary \ref{RN_Cor}, for each \ $\vartheta \in \RR$, \ $\alpha \in \RR$ \ and
 \ $T \in \RR_{++}$, \ we have
 \[
   \log\left(\frac{\dd\PP_{\vartheta+\alpha\hspace*{0.1mm}r_{\vartheta,T},\,T}}
                  {\dd\PP_{\vartheta,T}}
             (X^{(\vartheta)}|_{[-r,T]})\right)
   = \alpha \Delta_{\vartheta,T} - \frac{1}{2} \alpha^2 J_{\vartheta,T} .
 \]
Consequently, for each \ $(u, v) \in \RR^2$, \ we obtain
 \begin{align*}
  \EE\bigl(\exp\bigl\{\ii u \Delta_{\vartheta+\alpha\hspace*{0.1mm}r_{\vartheta,T},\,T}
                       + \ii v J_{\vartheta+\alpha\hspace*{0.1mm}r_{\vartheta,T},\,T}
                \bigr\}\bigr)
  &= \EE\biggl(\exp\bigl\{\ii u \Delta_{\vartheta,T} + \ii v J_{\vartheta,T}\bigr\}
               \frac{\dd\PP_{\vartheta+\alpha\hspace*{0.1mm}r_{\vartheta,T},\,T}}
                    {\dd\PP_{\vartheta,T}}
                (X^{(\vartheta)}|_{[-r,T]})\biggr) \\
  &=\EE\Bigl(\exp\Bigl\{(\ii u + \alpha) \Delta_{\vartheta,T}
                          + \Bigl(\ii v - \frac{1}{2} \alpha^2\Bigr)
                            J_{\vartheta,T}\Bigr\}\Bigr) .
 \end{align*}
By the continuous mapping theorem, \eqref{LAQD} yields
 \[
   \exp\Bigl\{(\ii u + \alpha) \Delta_{\vartheta,T}
              + \Bigl(\ii v - \frac{1}{2} \alpha^2\Bigr) J_{\vartheta,T}\Bigr\}
   \distr
   \exp\Bigl\{(\ii u + \alpha) \Delta_\vartheta
              + \Bigl(\ii v - \frac{1}{2} \alpha^2\Bigr) J_\vartheta\Bigr\}
 \]
 as \ $T \to \infty$.
\ The family
 \ $\bigl\{\exp\bigl\{(\ii u + \alpha) \Delta_{\vartheta,T}
                      + \bigl(\ii v - \frac{1}{2} \alpha^2\bigr) J_{\vartheta,T}\bigr\}
           : T \in \RR_+\bigr\}$
 \ is uniformly integrable, since for each \ $T \in \RR_+$, \ we have
 \[
   \Bigl|\exp\Bigl\{(\ii u + \alpha) \Delta_{\vartheta,T}
                    + \Bigl(\ii v - \frac{1}{2} \alpha^2\Bigr)
                      J_{\vartheta,T}\Bigr\}\Bigr|
   = \exp\Bigl\{\alpha \Delta_{\vartheta,T}
                - \frac{1}{2} \alpha^2 J_{\vartheta,T}\Bigr\} ,
 \]
 and
 \begin{equation}\label{mart}
   \EE\Bigl(\exp\Bigl\{\alpha \Delta_{\vartheta,T} 
                       - \frac{1}{2} \alpha^2 J_{\vartheta,T}\Bigr\}\Bigr)
   = \EE\biggl(\frac{\dd\PP_{\vartheta+\alpha\hspace*{0.1mm}r_{\vartheta,T},\,T}}
                    {\dd\PP_{\vartheta,T}}
               (X^{(\vartheta)}|_{[-r,T]})\biggr)
   = 1 ,
 \end{equation}
 since the process
 \ $\Bigl(\frac{\dd\PP_{\vartheta+\alpha\hspace*{0.1mm}r_{\vartheta,T},\,T}}
               {\dd\PP_{\vartheta,T}}
          (X^{(\vartheta)}|_{[-r,T]})\Bigr)_{T\in\RR_+}$
 \ is a martingale, see Lemma \ref{RN}.
Thus, by the moment convergence theorem,
 \begin{equation}\label{MCT}
  \begin{aligned}
   \EE\bigl(\exp\bigl\{\ii u
                       \Delta_{\vartheta+\alpha\hspace*{0.1mm}r_{\vartheta,T},\,T}
                       + \ii v J_{\vartheta+\alpha\hspace*{0.1mm}r_{\vartheta,T},\,T}
                \bigr\}\bigr)  
   &= \EE\Bigl(\exp\Bigl\{(\ii u + \alpha) \Delta_{\vartheta,T}
                          + \Bigl(\ii v - \frac{1}{2} \alpha^2\Bigr)
                            J_{\vartheta,T}\Bigr\}\Bigr) \\
   &\to
    \EE\Bigl(\exp\Bigl\{(\ii u + \alpha) \Delta_\vartheta
                        + \Bigl(\ii v - \frac{1}{2} \alpha^2\Bigr)
                          J_\vartheta\Bigr\}\Bigr)
  \end{aligned}
 \end{equation}
 as \ $T \to \infty$.
\ On the other hand, for each \ $\vartheta \in \RR$ \ and \ $\alpha \in \RR$, \ we have
 \[
   \log\left(\frac{\dd\PP^{(\alpha)}_\vartheta}{\dd\PP^{(0)}_\vartheta}
             (\bcX^{(0)}_\vartheta)\right)
   = \alpha \Delta^{(0)}_\vartheta - \frac{1}{2} \alpha^2 J^{(0)}_\vartheta .
 \]
Consequently, for each \ $(u, v) \in \RR^2$, \ we obtain
 \begin{align*}
  \EE\bigl(\exp\bigl\{\ii u \Delta^{(\alpha)}_\vartheta
                      + \ii v J^{(\alpha)}_\vartheta\bigr\}\bigr)
  &= \EE\biggl(\exp\bigl\{\ii u \Delta^{(0)}_\vartheta + \ii v J^{(0)}_\vartheta\bigr\}
               \frac{\dd\PP^{(\alpha)}_\vartheta}{\dd\PP^{(0)}_\vartheta}
                (\bcX^{(0)}_\vartheta)\biggr) \\
  &=\EE\Bigl(\exp\Bigl\{(\ii u + \alpha) \Delta^{(0)}_\vartheta
                        + \Bigl(\ii v - \frac{1}{2} \alpha^2\Bigr)
                          J^{(0)}_\vartheta\Bigr\}\Bigr) .
 \end{align*}
The aim of the following discussion is to show
 \ $(\Delta^{(0)}_\vartheta, J^{(0)}_\vartheta) = (\Delta_\vartheta, J_\vartheta)$.
\ For each \ $\lambda \in \Lambda_\vartheta \cap (\ii \RR)$, \eqref{SDElc} implies
 \ $\cX_{\vartheta,\lambda,0}^{(0)}(t) = \cW_{\Im(\lambda)}(t)$ \ and
 \ $\cX_{\vartheta,\lambda,\ell}^{(0)}(t)
    = \int_0^t \cX_{\vartheta,\lambda,\ell-1}^{(0)}(u) \, \dd u$
 \ for all \ $\ell \in \{1, \ldots, m^*_\vartheta\}$ \ and \ $t \in [0, 1]$.
\ By induction,
 \begin{equation}\label{Ito}
  \cX_{\vartheta,\lambda,\ell}^{(0)}(t)
  = \frac{1}{\ell!} \int_0^t (t - u)^\ell \, \dd\cW_{\Im(\lambda)}(u)
  = \cW_{\Im(\lambda),\ell}(t) , \qquad
  t \in [0, 1] ,
 \end{equation}
 for all \ $\ell \in \{1, \ldots, m^*_\vartheta\}$.
\ Indeed, by It\^o's formula,
 \[
   \cX_{\vartheta,\lambda,1}^{(0)}(t)
   = \int_0^t \cX_{\vartheta,\lambda,0}^{(0)}(u) \, \dd u
   = \int_0^t \cW_{\Im(\lambda)}(u) \, \dd u
   = \int_0^t (t - u) \, \dd \cW_{\Im(\lambda)}(u) ,
 \]
 hence we obtain \eqref{Ito} for \ $\ell = 1$.
\ If \eqref{Ito} holds for \ $\ell - 1$, \ then, again by It\^o's formula,
 \begin{align*}
  \cX_{\vartheta,\lambda,\ell}^{(0)}(t)
  &= \int_0^t \cX_{\vartheta,\lambda,\ell-1}^{(0)}(s) \, \dd s
   = \int_0^t
      \left(\frac{1}{(\ell-1)!}
            \int_0^s (s - u)^{\ell-1} \, \dd\cW_{\Im(\lambda)}(u)\right)
      \dd s \\
  &= \int_0^t
      \left(\frac{1}{(\ell-2)!}
            \int_0^s (s - u)^{\ell-2} \cW_{\Im(\lambda)}(u) \, \dd u\right)
      \dd s \\
  &= \int_0^t
      \left(\frac{1}{(\ell-2)!} \int_0^s (s - u)^{\ell-2} \dd s\right)
      \cW_{\Im(\lambda)}(u) \, \dd u \\
  &= \int_0^t \frac{1}{(\ell-1)!} (t - u)^{\ell-1} \cW_{\Im(\lambda)}(u) \, \dd u
   = \frac{1}{\ell!} \int_0^t (t - u)^\ell \, \dd\cW_{\Im(\lambda)}(u) ,
 \end{align*}
 thus we obtain \eqref{Ito} for \ $\ell$, \ and we conclude
 \begin{gather*}
  \Delta^{(0)}_{\vartheta,\lambda}
  = \begin{cases}
     c_{\vartheta,0,m^*_\vartheta} \int_0^1 \cW_{0,m^*_\vartheta}(t) \, \dd \cW_0(t)
      & \text{if \ $\Im(\lambda) = 0$,} \\
     2 \Re\bigl(\int_0^1
                 c_{\vartheta,\lambda,m^*_\vartheta}
                 \cW_{\Im(\lambda),m^*_\vartheta}(t)
                 \, \dd \overline{\cW_{\Im(\lambda)}(t)}\bigr)
      & \text{if \ $\Im(\lambda) \in \RR_{++}$,}
    \end{cases} \\
  J^{(0)}_{\vartheta,\lambda}
  = \begin{cases}
     c_{\vartheta,0,m^*_\vartheta}^2 \int_0^1 \cW_{0,m^*_\vartheta}(t)^2 \, \dd t
      & \text{if \ $\Im(\lambda) = 0$,} \\
     2 |c_{\vartheta,\lambda,m^*_\vartheta}|^2
     \int_0^1 |\cW_{\Im(\lambda),m^*_\vartheta}(t)|^2 \, \dd t
      & \text{if \ $\Im(\lambda) \in \RR_{++}$.}
    \end{cases}
 \end{gather*}
Using the identity \eqref{conj}, we obtain
 \ $(\Delta^{(0)}_\vartheta, J^{(0)}_\vartheta) = (\Delta_\vartheta, J_\vartheta)$,
 \ and hence, \eqref{MCT} implies
 \begin{align*}
  &\EE\bigl(\exp\bigl\{\ii u
                       \Delta_{\vartheta+\alpha\hspace*{0.1mm}r_{\vartheta,T},\,T}
                       + \ii v J_{\vartheta+\alpha\hspace*{0.1mm}r_{\vartheta,T},\,T}
                \bigr\}\bigr)
   \to \EE\Bigl(\exp\Bigl\{(\ii u + \alpha) \Delta_\vartheta
                           + \Bigl(\ii v - \frac{1}{2} \alpha^2\Bigr)
                             J_\vartheta\Bigr\}\Bigr) \\
  &= \EE\Bigl(\exp\Bigl\{(\ii u + \alpha) \Delta^{(0)}_\vartheta
                           + \Bigl(\ii v - \frac{1}{2} \alpha^2\Bigr)
                             J^{(0)}_\vartheta\Bigr\}\Bigr)
   = \EE\bigl(\exp\bigl\{\ii u \Delta^{(\alpha)}_\vartheta
                         + \ii v J^{(\alpha)}_\vartheta\bigr\}\bigr)
 \end{align*}
 as \ $T \to \infty$, \ and the continuity theorem yields \eqref{convDJ}.
 
The rest of the statement can be proved as Theorem 2.10 in our paper
 \cite{BenPap2017b}. 
By Corollary \ref{RN_Cor}, convergence \eqref{convDJ} and the continuous mapping
 theorem, for each \ $\vartheta \in \RR$, \ $\alpha \in \RR$ \ and \ $h \in \RR$, \ we
 obtain
 \begin{align*}
  \frac{\dd\PP_{\vartheta+(\alpha+h)r_{\vartheta,T},\,T}}
       {\dd\PP_{\vartheta+\alpha\hspace*{0.1mm}r_{\vartheta,T},\,T}}
   (X^{(\vartheta+\alpha\hspace*{0.1mm}r_{\vartheta,T})}|_{[-r,T]})
  &= \exp\Bigl\{h \Delta_{\vartheta+\alpha\hspace*{0.1mm}r_{\vartheta,T},\,T}
                - \frac{1}{2} h^2
                  J_{\vartheta+\alpha\hspace*{0.1mm}r_{\vartheta,T},\,T}\Bigr\} \\
  &\distr 
   \exp\Bigl\{h \Delta^{(\alpha)}_\vartheta
              - \frac{1}{2} h^2 J^{(\alpha)}_\vartheta\Bigr\} ,
   \qquad \text{as \ $T \to \infty$.}
 \end{align*}
We have
 \ $\EE\bigl(\exp\bigl\{h \Delta^{(\alpha)}_\vartheta
                        - \frac{1}{2} h^2 J^{(\alpha)}_\vartheta\bigr\}\bigr)
    = \EE\Bigl(\frac{\dd\PP^{(\alpha+h)}_\theta}{\dd\PP^{(\alpha)}_\vartheta}
                (\bcX^{(\alpha)})\Bigr)
    = 1$,
 \ as in \eqref{mart}.
By Le Cam's first lemma, we conclude that the families
 \ $(\PP_{\vartheta+(\alpha+h)r_{\vartheta,T},\,T})_{T\in\RR_+}$ \ and
 \ $(\PP_{\vartheta+\alpha\hspace*{0.1mm}r_{\vartheta,T},\,T})_{T\in\RR_+}$ \ are
 mutually contiguous as \ $T \to \infty$.
\ Therefore, for each \ $h, h_0 \in \RR$, \ the probability of the set on which we have
 \begin{align*}
  &\log\left(\frac{\dd\PP_{\vartheta+(\alpha+h)r_{\vartheta,T},\,T}}
                  {\dd\PP_{\vartheta+(\alpha+h_0)r_{\vartheta,T},\,T}}
             (X^{(\vartheta+(\alpha+h_0)r_{\vartheta,T})}|_{[-r,T]})\right) \\
  &= \log\left(\frac{\dd\PP_{\vartheta+(\alpha+h)r_{\vartheta,T},\,T}}
                    {\dd\PP_{\vartheta+\alpha\hspace*{0.1mm}r_{\vartheta,T},\,T}}
                (X^{(\vartheta+\alpha\hspace*{0.1mm}r_{\vartheta,T})}|_{[-r,T]})\right)
     - \log\left(\frac{\dd\PP_{\vartheta+(\alpha+h_0)r_{\vartheta,T},\,T}}
                      {\dd\PP_{\vartheta+\alpha\hspace*{0.1mm}r_{\vartheta,T},\,T}}
                (X^{(\vartheta+\alpha\hspace*{0.1mm}r_{\vartheta,T})}|_{[-r,T]})\right)
 \end{align*}
 converges to one as \ $T \to \infty$.
\ Hence, by Lemma \ref{RN}, we have
 \[
   \log\left(\frac{\dd\PP_{\vartheta+(\alpha+h)r_{\vartheta,T},\,T}}
                  {\dd\PP_{\vartheta+(\alpha+h_0)r_{\vartheta,T},\,T}}
             (X^{(\vartheta+(\alpha+h_0)r_{\vartheta,T})}|_{[-r,T]})\right)
   = (h - h_0) \Delta_{\vartheta+\alpha\hspace*{0.1mm}r_{\vartheta,T},\,T}
     - \frac{1}{2} (h^2 - h_0^2)
       J_{\vartheta+\alpha\hspace*{0.1mm}r_{\vartheta,T},\,T} 
 \]
 with probability converging to one as \ $T \to \infty$.
\ By Le Cam's third lemma, for each \ $h \in \RR$, \ the distribution of
 \ $(\Delta_{\vartheta+\alpha\hspace*{0.1mm}r_{\vartheta,T},\,T},
     J_{\vartheta+\alpha\hspace*{0.1mm}r_{\vartheta,T},\,T})$
 \ under \ $\PP_{\vartheta+(\alpha+h)r_{\vartheta,T},\,T}$
 \ tends to \ $\QQ^{(\alpha)}_{\vartheta,h}$ \ as \ $T \to \infty$, \ hence we obtain
 the statement.    
\proofend

\section{Asymptotics of the maximum likelihood estimators}
\label{AMLE}

For fixed \ $T \in \RR_{++}$, \ maximizing
 \ $\log\bigl(\frac{\dd\PP_{\vartheta,T}}{\dd\PP_{\vartheta_0,T}}
               (X^{(\vartheta_0)}|_{[-r,T]})\bigr)$
 \ in \ $\vartheta \in \RR$ \ and then replacing \ $X^{(\vartheta_0)}$ \ by
 \ $X^{(\vartheta)}$ \ gives the MLE of \ $\vartheta$ \ based on the observations
 \ $(X^{(\vartheta)}(t))_{t\in[-r,T]}$ \ having the form
 \[
   \hvartheta_T
   = \frac{\int_0^T Y^{(\vartheta)}(t) \, \dd X^{(\vartheta)}(t)}
          {\int_0^T Y^{(\vartheta)}(t)^2 \, \dd t} ,
 \]
 provided that \ $\int_0^T Y^{(\vartheta)}(t)^2 \, \dd t > 0$.
\ Using the SDDE \eqref{SDDE}, one can check that
 \[
   \hvartheta_T - \vartheta
   = \frac{\int_0^T Y^{(\vartheta)}(t) \, \dd W(t)}
          {\int_0^T Y^{(\vartheta)}(t)^2 \, \dd t} ,
 \]
 hence
 \[
   r_{\vartheta,T}^{-1} (\hvartheta_T - \vartheta)
   = \frac{\Delta_{\vartheta,T}}{J_{\vartheta,T}} .
 \]
If \ $v_\vartheta^* = 0$, \ then \eqref{convDJ} and the Continuous Mapping Theorem
 yield
 \[
   r_{\vartheta,T}^{-1} (\hvartheta_T - \vartheta)
   \distr \frac{\Delta_\vartheta}{J_\vartheta} \qquad \text{as \ $T \to \infty$,}
 \]
 since \ $\PP(J_\vartheta > 0) = 1$.

The MLE of the parameter \ $\alpha$ \ based on the observations
 \ $(X^{(\vartheta+\alpha\hspace*{0.1mm}r_{\vartheta,T})}(t))_{t\in[-r,T]}$ \ has the
 form
 \ $\halpha_T
    = r_{\vartheta,T}^{-1}
      ((\vartheta + \alpha\hspace*{0.1mm}r_{\vartheta,T})\hspace*{1mm}\widehat{}
       - \vartheta)$,
 \ where \ $(\vartheta + \alpha\hspace*{0.1mm}r_{\vartheta,T})\hspace*{1mm}\widehat{}$
 \ denotes the MLE of \ $\vartheta + \alpha\hspace*{0.1mm}r_{\vartheta,T}$.
\ Hence
 \[
   \halpha_T - \alpha
   = r_{\vartheta,T}^{-1}
     ((\vartheta + \alpha\hspace*{0.1mm}r_{\vartheta,T})\hspace*{1mm}\widehat{}
      - \vartheta)
     - \alpha
   = r_{\vartheta,T}^{-1}
     ((\vartheta + \alpha\hspace*{0.1mm}r_{\vartheta,T})\hspace*{1mm}\widehat{}
      - (\vartheta + \alpha\hspace*{0.1mm}r_{\vartheta,T}))
   = \frac{\Delta_{\vartheta + \alpha\hspace*{0.1mm}r_{\vartheta,T}}}
          {J_{\vartheta + \alpha\hspace*{0.1mm}r_{\vartheta,T}}} .
 \]
Maximizing
 \ $\log\Bigl(\frac{\dd\PP^{(\alpha)}_\vartheta}{\dd\PP^{(\alpha_0)}_\vartheta}
               \bigl(\bcX^{(\alpha_0)}_\vartheta\bigr)\Bigr)$
 \ in \ $\alpha \in \RR$ \ and then replacing \ $\bcX^{(\alpha_0)}_\vartheta$ \ by
 \ $\bcX^{(\alpha)}_\vartheta$ \ gives the MLE of \ $\alpha$ \ based on the
 observations \ $(\bcX^{(\alpha)}_\vartheta(t))_{t\in[0,1]}$ \ having the form
 \[
   \halpha
   = \frac{\sum\limits_{\underset{\SC\tm_\vartheta(\lambda)=m_\vartheta^*}
                                 {\lambda \in \Lambda_\vartheta\cap(\ii\RR)}}
            c_{\vartheta,\lambda,m_\vartheta^*}
            \int_0^1
             \cX^{(\alpha)}_{\vartheta,\lambda,m_\vartheta^*}(s)
             \, \dd \overline{\cX^{(\alpha)}_{\vartheta,\lambda,0}(s)}}
          {\sum\limits_{\underset{\SC\tm_\vartheta(\lambda)=m_\vartheta^*}
                                 {\lambda \in \Lambda_\vartheta\cap(\ii\RR)}}
            |c_{\vartheta,\lambda,m_\vartheta^*}|^2
            \int_0^1
             |\cX^{(\alpha)}_{\vartheta,\lambda,m_\vartheta^*}(s)|^2 \, \dd s} ,
 \]
 provided that the denominator is not zero, which is satisfied with probability one,
 since the coefficient \ $c_{\vartheta,\lambda,m_\vartheta^*}$ \ is not zero, and
 \ $(\cX^{(\alpha)}_{\vartheta,\lambda,m_\vartheta^*}(s))_{s\in[0,1]}$ \ is a
 non-degenerate Gaussian process.
Using the SDE \eqref{SDElc}, one can check that
 \[
   \halpha - \alpha
   = \frac{\sum\limits_{\underset{\SC\tm_\vartheta(\lambda)=m_\vartheta^*}
                                 {\lambda \in \Lambda_\vartheta\cap(\ii\RR)}}
            c_{\vartheta,\lambda,m_\vartheta^*}
            \int_0^1
             \cX^{(\alpha)}_{\vartheta,\lambda,m_\vartheta^*}(s)
             \, \dd \overline{\cW_{\Im(\lambda)}(s)}}
          {\sum\limits_{\underset{\SC\tm_\vartheta(\lambda)=m_\vartheta^*}
                                 {\lambda \in \Lambda_\vartheta\cap(\ii\RR)}}
            |c_{\vartheta,\lambda,m_\vartheta^*}|^2
            \int_0^1
             |\cX^{(\alpha)}_{\vartheta,\lambda,m_\vartheta^*}(s)|^2 \, \dd s}
   = \frac{\Delta^{(\alpha)}_\vartheta}{J^{(\alpha)}_\vartheta} .
 \]
Thus \eqref{convDJ} and the continuous mapping theorem yield the following result.
 
\begin{Thm}\label{MLEconv}
Let \ $\vartheta \in \RR$ \ with \ $v_\vartheta^* = 0$.
\ For each \ $T \in \RR_{++}$, \ let \ $\halpha_T$ \ denote the MLE of \ $\alpha$
 \ based on the observations
 \ $(X^{(\vartheta+\alpha\hspace*{0.1mm}r_{\vartheta,T})}(t))_{t\in[-r,T]}$.
\ Let \ $\halpha$ \ denote the MLE of \ $\alpha$ \ based on the observations
 \ $(\bcX^{(\alpha)}_\vartheta(t))_{t\in[0,1]}$.
\ Then \ $\halpha_T \distr \halpha$ \ as \ $T \to \infty$.
\end{Thm}


\begin{thebibliography}{9}

\bibitem{And}
\textsc{Anderson, T. W.} (1959).
On asymptotic distributions of estimates of parameters of stochastic difference equations.
\textit{Ann.\ Math.\ Statist.}
\textbf{30(3)} 676--687.

\bibitem{Arato}
\textsc{Arat\'o, M.} (1982).
\textit{Linear Stochastic Systems with Constant Coefficients. A Statistical Approach.}
Lectures Notes in Control and Information Sciences, vol. 45.
Berlin: Springer-Verlag. (In Russian, Moscow: Nauka, 1989.)

\bibitem{BenPap2017a}
\textsc{Benke, J. M.} and \textsc{Pap, G.} (2017).
One-parameter statistical model for linear stochastic differential equation
 with time delay.
\textit{Statistics}
\textbf{51(3)} 510--531.

\bibitem{BenPap2017b}
\textsc{Benke, J. M.} and \textsc{Pap, G.} (2017).
Local asymptotic quadraticity of statistical experiments connected with a Heston model.
\textit{Acta Sci.\ Math.\ (Szeged)}
\textbf{83(1-2)} 313--344.

\bibitem{Bob}
\textsc{Bobkoski, J. M.} (1983).
Hypothesis testing in nonstationary time series.
PhD thesis, University of Wisconsin, Madison.

\bibitem{BucCha}
\textsc{Buchmann, B.} and \textsc{Chan, N. H.} (2013).
Unified asymptotic theory for nearly unstable AR($p$) processes.
\textit{Stochastic Process.\ Appl.} 
\textbf{123} 952--985.

\bibitem{Chan}
\textsc{Chan, N. H.} (2009).
Time Series with Roots on or Near the Unit Circle. \text{In:} Mikosch ,T., Krei{\ss}, J. P., Davis, R., Andersen, T. (Eds): \textit{Handbook of Financial Time Series}, 695--707. Springer, Berlin, Heidelberg.

\bibitem{ChanWei}
\textsc{Chan, N. H.} and \textsc{Wei, C. Z.} (1987).
Asymptotic inference for nearly nonstationary AR(1) processes.
\textit{Ann.\ Statist.} 
\textbf{15(3)} 1050--1063.

\bibitem{ChanWei2}
\textsc{Chan, N. H.} and \textsc{Wei, C. Z.} (1988).
The parameter inference for nearly nonstationary time series.
\textit{J. Amer.\ Statist.\ Assoc.} 
\textbf{83} 857--862.

\bibitem{GreWef}
\textsc{Greenwood, P. E.} and \textsc{Wefelmeyer, W.} (1993).
Asymptotic minimax results for stochastic process families with critical points.
\textit{Stochastic Process.\ Appl.} 
\textbf{44} 107--116.

\bibitem{JSh}
\textsc{Jacod, J.} and \textsc{Shiryaev, A. N.} (2003).
\textit{Limit Theorems for Stochastic Processes}, 2nd ed.
Springer-Verlag, Berlin.

\bibitem{Jeg}
\textsc{Jeganathan, P.} (1991).
On the asymptotic behavior of least-squares estimators in AR time series with
roots near the unit circle.
\textit{Econometric Theory}
\textbf{7} 269--306.

\bibitem{LipShiI}
\textsc{Liptser, R. S.} and \textsc{Shiryaev, A. N.} (2001).
\textit{Statistics of Random Processes I. General theory}, 2nd edition.
Springer-Verlag, Berlin, Heidelberg.

\bibitem{ManWal}
\textsc{Mann, H. B.} and \textsc{Wald, A.} (1943).
On the statistical treatment of linear stochastic difference equations.
\textit{Econometrica}
\textbf{11(3)} 173--220.

\bibitem{MeerPapZui}
\textsc{van der Meer, T.}, \textsc{Pap, G.} and \textsc{van Zuijlen, M. C. A.} (1999).
Asymptotic inference for nearly unstabel AR(p) processes.
\textit{Econometric Theory}
\textbf{15(2)} 184--217.

\bibitem{Phi}
\textsc{Phillips, P. C. B.} (1987).
Towards a unified asymptotic theory for autoregression.
\textit{Biometrika}
\textbf{74(3)} 535--547.

\bibitem{Phi1}
\textsc{Phillips, P. C. B.} (1987).
Time series regression with a unit root.
\textit{Econometrica}
\textbf{55} 277--301.

\bibitem{Phi2}
\textsc{Phillips, P. C. B.} (1989).
Partially identified econometric models.
\textit{Econometric Theory}
\textbf{5} 181--240.

\bibitem{White}
\textsc{White, J. S.} (1958).
The limiting distribution of the serial correlation in the explosive case.
\textit{Ann.\ Math.\ Statist.}
\textbf{29(4)} 1188--1197.

\end{thebibliography}
\end{document}